\documentclass[12pt]{amsart}
\usepackage{epsfig,amssymb,epic,eepic,color}
\usepackage[all]{xy}
 \usepackage[T1]{fontenc}        
\newtheorem{thm}{Theorem}[section]

\newtheorem{prop}[thm]{Proposition}
\newtheorem{defn}[thm]{Definition}
\newtheorem{lemme}[thm]{Lemma}
\newtheorem{cor}[thm]{Corollary}

\newtheorem{remarque}[thm]{Remark}

\newtheorem{rien}[thm]{}
 
\newcommand{\be}{\begin{enumerate}}
\newcommand{\ee}{\end{enumerate}}
\newcommand{\bi}{\begin{itemize}}
\newcommand{\ei}{\end{itemize}}

\def\R{\mathbb{R}}

\def\Z{\mathbb{Z}}

\def\F{\mathbb{F}}

\def\a{\alpha}

\def\ga{\gamma}

\def\al{\alpha}
\def\be{\beta}
 
\def\de{\delta}
\def\De{\Delta}

\def\la{\lambda}

\def\si{\sigma}
\def\Si{\Sigma}

\def\ep{\varepsilon}

\def\ds{\displaystyle}

\def\nd{\noindent}
\def\bull{\hfill$\Box$\\}
\def\proof{\nd {\bf Proof.\ }}

\begin{document}
\vskip 1cm
\begin{center}
{\sc On an article by S. A. Barannikov}


\vspace{1cm}
{\sc Fran\c cois Laudenbach}
\end{center}

\title{}
\author{}
\address{Laboratoire de
math\'ematiques Jean Leray,  UMR 6629 du
CNRS, Facult\'e des Sciences et Techniques,
Universit\'e de Nantes, 2, rue de la
Houssini\`ere, F-44322 Nantes cedex 3,
France.}
\email{francois.laudenbach@univ-nantes.fr}



\begin{abstract} Given a Morse function $f$ on a closed manifold $M$ with distinct critical values,
and given a field $\F$, there is a canonical complex, called the Morse-Barannikov complex,
which is  {\it equivalent} to any Morse complex 
associated with $f$ and whose form is {\it simple}. In particular, the homology of $M$ with coefficients in $\F$ is immediately readable on this complex. The bifurcation theory of 
this complex in a generic one-parameter family of functions will be investigated. Applications 
to the boundary manifolds will be given.
\end{abstract}
\maketitle

\thispagestyle{empty}
\vskip 1cm
\medskip
Here is an expanded version of the lectures given in the Winter school organized 
in La Llagonne 
by the University Paul Sabatier in Toulouse (January 2013)\footnote{This text was written in February 2013. It replaces the previous version posted by error.}.
There are three parts:
\begin{itemize}

\item[I)] The Morse-Barannikov complex (after C. Viterbo),

\item[II)] Bifurcations, 

\item[III)] The non-empty boundary case. 
\end{itemize}

 \vskip 1cm
 \section{The Morse-Barannikov complex}
 
 We adopt a presentation which is a mix of the presentation given by S. Barannikov
 in \cite{bara} and a more abstract one given by C. Viterbo in
 his joint work with D. Le Peutrec and F. Nier \cite{viterbo} which we
 slightly simplify.
 
 We are given a closed manifold $M$, a Morse function $f: M\to \R$ whose critical values are
 distinct, and a field $\F$. For each integer $k$ the critical points of index $k$ 
 are numbered  in the increasing order
 of the critical values: $f(p_1)<f(p_2)<\ldots$ (the function is just generic and it is not assumed to be ordered). We shall often identify the set of critical points and the set of critical values.
 
 For defining the {\it Morse complex} it is necessary to have two extra data:
 \begin{itemize}
 \item A (decreasing) {\it pseudo-gradient}, that is, a vector field on $M$ which satifies
 $X\cdot f<0$ out of the critical points, and some non-degeneracy condition for the vanishing 
 of $X$ at each critical point; therefore, the zeroes of $X$ are hyperbolic. As a consequence,
 each critical point $p$ has a stable manifold $W^s(p)$ and an unstable manifold $W^u(p)$. 
 This pseudo-gradient is chosen Morse-Smale, 
  a generic property meaning that the stable manifolds are transverse to the unstable manifolds.
 
 \item An orientation of the unstable manifolds.
  \end{itemize}
  
  The Morse complex $C_*(f,X)$ is made as follows. In degree $k$ the module $C_k(f,X)$
  is the free $\Z$-module generated by the critical points of index $k$ and the differential
  $\partial_{k}: C_{k}(f,X)\to C_{k-1}(f,X)$ counts the signed number of connecting orbits.
  Observe that, for $(p, q)$, a pair of critical points of respective index $k$ and $k-1$,
  $W^u(p)\cap W^s(q)$ is made of a finite number of connecting orbits. Since $W^s(q)$ is co-oriented by the orientation of 
  $W^u(q)$, each orbit descending from $p$ to $q$ gets a sign. Define $a_p^q$ to be the signed 
  number of the connecting orbits and define the Morse differential $C_k(f,X)$ by 
 $$\partial_k(p)=\sum a_p^q q\,.
 $$
 
\begin{thm}{\bf (Milnor \cite{h-cob}, Th. 7.2)}. The Morse complex is a chain complex: $\partial\circ \partial=0$.
 Moreover, $H\left(C_*(f,X)\right)\cong H_*(M;\Z)$. A fortiori, 
 $H\left(C_*(f,X); \F\right)\cong H_*(M;\F)$.
 \end{thm}
 
\begin{defn} A chain complex  $C_*$ with coeffcients in $\F$
is said to be \emph{$\F$-equivalent} to $C_*(f,X)$ if it has the same generators 
and if its differential $\delta$ is made from $\partial $ by conjugating in each degree
by an invertible upper triangular matrix $T$ with coefficients in $\F$:
\[\xymatrix{&C_{k+1}(f,X)\otimes\F\ar[d]_{T} \ar[r]^\partial& C_k (f,X)\otimes\F\ar[r]^{\partial}
\ar[d]_{T}&
C_{k-1}(f,X)\otimes\F\ar[d]_T\\
&C_{k+1}\ar[r]^{\de} &C_k\ar[r]^{\de}&C_{k-1}
}\]\\
Here, $C_k (f,X)\otimes\F$ is a vector space equipped with its canonical ordered basis.
\end{defn}
 
 \begin{remarque} {\rm When changing the pseudo-gradient or the orientations of the unstable manifolds
 the Morse complex is changed by $\Z$-equivalence. Conversely, if $\dim M>1$ and 
 if the level sets of $f$ are connected (or $f$ has one local minimum 
 and one local maximum only, {\it i.e.} $f$ is {\it polar} in Morse's terminology) every $\Z$-equivalence is realizable
 by such changes. When the  coefficients are in a field an $\F$-equivalence has no longer 
 such a geometrical meaning in general. But, an $\F$-equivalence keeps the memory
  of the filtration by the sub-level sets of the function $f$. This fact will be used in the last step 
  of the proof of Barannikov's theorem.
  }
 \end{remarque}
 
 \begin{thm} {\bf (Barannikov \cite{bara})} \label{main} The Morse complex $C_*(f,X)$ is 
 $\F$-equivalent to a \emph{simple}
 complex $\left(C_*,\partial_B\right)$, that is, for every generator $p$, 
 $\partial_B(p)$ is 0 or a generator and $\partial_B(p)\not=\partial_B(p')$ if $p\not=p'$ and $\partial_B(p)\not=0$. Moreover, $\left(C_*,\partial_B\right)$ is unique and depends only
 on $C_*(f,X)\otimes \F$ for any pseudo-gradient $X$.  \end{thm}
 
 This complex  is called the {\it Morse-Barannikov complex} associated with the Morse function $f$; 
  it depends on the field $\F$.

 \begin{cor} The homology $H_*(M;\F)$ is graded isomorphic to the sub-space generated by the critical points having the homological type, in the sense given below: $\partial_B (p)=0$ and
 $p\notin {\rm Im}\,\partial_B$.  
 \end{cor}
 
 One  important point in the statement is the coupling of some  critical points, the unpaired generators being ``isolated'' in the complex. This fact plays a deep r\^ole in the work by Y. Chekanov
 \& P. Pushkar \cite{chekanov}. When the Morse complex is concentrated in two degrees,
 the statement amounts to the fact that the double coset $GL(n,\Z)/T(n)\times T(n)$ is isomophic
 to the symmetric group $S_n$, a fact which was important in Cerf's work on pseudo-isotopy (here, $T(n)$ denotes the sub-group of invertible upper triangular matrices) \cite{cerf}.
 
 In Barannikov's paper the proof  of existence follows more or less from Gauss' algorithm. The proof of uniqueness remains mysterious. It is clarified by C. Viterbo in \cite{viterbo}.

 Viterbo's important remark is that the critical points of the given Morse function $f$
 are divided in three types: {\it upper, lower} and {\it homological}, depending of the place of 
  a zero map 
 in the diagram below of $\F$-vector spaces and 
 $\F$-linear maps, in which $c$ denotes a critical value of index $k+1$ and, for brevity, 
 $c+\ep$ stands for the sub-level set $f^{c+\ep}:= f^{-1}\bigl((-\infty, c+\ep]\bigr)$:
 \[\xymatrix{&&&\F&&\\
 0 \ar[r] &H_{k+1}(c-\ep) \ar[r]  &H_{k+1}(c+\ep) \ar[r]^-{J} &
  H_{k+1}(c+\ep,c-\ep) \ar[u]_{\cong}\ar[r]^-{\Delta} \ar[d]^I 
  &H_k(c-\ep)\ar[r] &H_k(c+\ep)\ar[r] &0\\
  &&& H_{k+1}(+\infty,c-\ep)\ar[ur]_{\De'}&&
}\]\\
The horizontal line is an exact sequence.
 The critical point $p$ such that $f(p)=c$ is said to be of {\it upper type} when $J=0$, implying $
 \De$ injective. It is said to be of {\it lower type} when $J$ is surjective and $I=0$.
  It is said of {\it homological type} when $I$ is injective and $\De=0$. Clearly, since $\F$ is a field, 
 all possibilities are covered, making a partition of the critical points.
 
 The type of a critical point $p$ of $f$ is readable on the Morse complex with coefficients
 in $\F$, that is $C_*(f,X)\otimes \F$. For instance, $p$ is of lower type  if there is an $\F$-linear
  combination  of  critical points higher than $p$ whose boundary is $p$\,; this property does
  not depend on the chosen pseudo-gradient $X$.\\
 
 \nd {\sc The coupling of critical points.} If $p$ is a critical point of index $k+1$, the local 
 unstable manifold $W^u_{loc}(p)$ is unique up to isotopy and orientation. 
 Set $[p]:=\left[ W^u_{loc}(p)\right]\in H_{k+1}(c+\ep,c-\ep)  $
 where $c=f(p)$.
 If $p$ is of  upper type, one defines 
 $$\lambda(p): = \inf_\si\max(f\vert\si)
 $$
 where $\si$ runs among the $k$-cycles of the sub-level  set $f^{c-\ep}$ representing 
 $\De ([p])$; it is a critical value $\la(p)=f(q)$ where $q$ is a critical point of index $k$. 
 By identifying critical point and critical value, we set $q:=\la(p)$.
 
 \begin{lemme} \label{lem1} The critical point $q:=\la(p)$ is of lower type.
 \end{lemme}
 \nd \proof
 Denote $\De_q([p])$ the class of $\De([p])$ modulo the sub-level set $f(q)-\ep$ in \break
 $H_k\bigl(f(p)-\ep, f(q)-\ep\bigr)$. By definition of the minimax, this class is not zero
 and we have 
 $$\De_q([p])=\al I_q^p\bigl([q]\bigr), \al\in\F,\al\not= 0,
 $$ where $I_q^p: H_k(f(q)+\ep, f(q)-\ep)\to H_k(f(p)-\ep, f(q)-\ep)$ is induced by the inclusion.
 Thus, $W^u_{loc}(q)$ is the boundary of $\frac 1\al W^u(p)$ in the pair $(+\infty, f(q)-\ep)$.
 Hence $I([q])=0$.\bull
 
\nd {\sc The Barannikov differential}. Set $\partial_B(p)=\la(p)$ if $p$ is of upper type
and $\partial_B(p)=0$ in the two other cases. According to the previous lemma, 
$\partial_B\circ \partial_B=0$.

\begin{lemme} \label{surj}The map $\la$ defines a bijection from the set of critical points
of upper type onto the set of critical points of lower type.
\end{lemme}
\nd \proof 1) {\sc Injectivity.} Let $q=\partial_B(p)=\partial_B(p')$
 with $f(p)>f(p')$.
We are using the same notation as in Lemma \ref{lem1}. We have 
$\De_q([p])=\al I_q^p([q])\not=0$ in $H_k\bigl(f(p)-\ep, f(q)-\ep\bigr)$ 
and $\De_q([p'])=\al' I_{q}^{p'}([q])\not=0$ in $H_k\bigl(f(p')-\ep, f(q)-\ep\bigr)$ .
 By construction 
$I_q^p\De_q([p'])=0$. Thus, we have: 
$$\begin{array}{lcl}
\De_q([p]) &=& \De_q([p])-\ds{\frac{\al}{\al'}} I_q^p\De_q([p'])\\
&{}&\\
 &=&\al I_q^p([q]) - \ds{\frac{\al}{\al'}\al'}I_q^{p }([q])=0\,.
\end{array}
$$
Therefore, $f(q)$ is not the minimax value associated with $p$.\\

\nd 2) {\sc Surjectivity.} Let $q$ be a critical point of lower type and
 index $k$.
Set
$$\mu(q)=\inf_\si\max (f\vert\si)
$$where $\si$ runs among the relative chains of $({+\infty},f(q)-\ep)$
 whose boundary
is a relative cycle representing the class 
$[q]$. This $\mu(q)$ is a critical value of index $k+1$ with $\mu(q)=f(p)$ for some critical point $p$.
A chain $\si $ approximating  the infimum has a non vanishing class in 
$H_{k+1}\bigl(f(p)+\ep,f(p)-\ep\bigr)$; hence,  $[\si]=\beta [p]$
in the pair $(f(p)+\ep,f(p)-\ep)$ with
 $\beta\in \F, \ \beta\not=0$.
We have $\ds{\De_q\bigl([p]\bigr)}= I_q^p(\frac 1\beta [q])$. If this element is zero, this means that there is another relative chain bounded by $W^u(q)$
under the level of $p$, contradicting the definition of $\mu(q)$.
 Then, $\De_q([p])\not=0$.
{\it A fortiori}, $\De([p])\not=0$ and $p$ is of upper type.

We have also to show that $\la(p)=q$. 
By the above $\si$, $\De([p])$ is homologous to $[q]$ up a non zero scalar.
If it is homologous to $[q']$ with $f(q')<f(q)$, then $I_q^p([q])=0$,
 and this is not the case. 
\bull

At this point we have the uniqueness part in Barannikov's theorem.

\begin{lemme} The Morse-Barannikov complex is $\F$-equivalent to the Morse complex. In particular,
its homology is isomorphic to $H_*(M,\F)$.
\end{lemme}

\nd \proof Suppose we have a chain complex $(C_*,\partial)$,   $\F$-equivalent to the Morse complex and 
which is simple until the degree $k$. Then, $\partial\left(C_{k+1}\right)$ is orthogonal
to the critical points of upper type (with respect to the canonical scalar product of a
based vector space). If not, $\partial\circ\partial\not= 0$.

Let $p_1,\ldots, p_m$ be the critical points of index $k+1$, with $f(p_1)<\ldots< f(p_m)$.
Let $q_1,\ldots, q_r$ be the critical points of index $k$ whose type is lower or homological;
$f(q_1)<\ldots< f(q_r)$. We assume that,  for some $j\leq r$, we have  $\partial p_i=0$ or 
$\partial p_i=q_{k(i)} $
for every $i<j$, the map  $i\mapsto k(i)$ being injective.
We have 
$$
\partial p_j=\sum_{i<j,\partial p_i\not=0
}
\al_i\, q_{k(i)}+ \left\{\begin{array}{l}
0\\
{}\\
{\rm or}\\
{}\\
\beta_{k_0}\,q_{k_0} +\mathop{ \sum}\limits_{k<k_0,k\not=k(i)
}
\beta_k\,q_k,\quad {\rm with }\quad \beta_{k_0}\not=0\,.
\end{array}\right.
$$

In the first case we use the following upper triangular matrix in degree $k+1$:
$$\begin{array}{l}
T (p_j)=p_j-\ds{\sum_{i<j,\partial p_i\not=0}}\al_i\, p_i\\
T(p_\ell)= p_\ell \quad {\rm if}\quad \ell\not=j
\end{array}
$$ and we set $\bar \partial= \partial \circ T$. We get $\bar\partial (p_j)=0$ and $\bar\partial(p_i)=
\partial(p_i)$ for every $i<j$;
so we have improved the simpleness of the differential.

In the second case, we use an upper triangular matrix in both degree $k+1$ and $k$:
$$\begin{array}{l}
T(p_j)= \frac {1}{\beta_{k_0}}\left( p_j-\ds{\sum_{i<j,\, \partial p_i\not= 0}}\a_i\,p_i\right),\\
T(p_\ell)= p_\ell \quad {\rm if}\quad \ell\not=j,\\
T(q_{k_0})= q_{k_0} + \ds{\mathop{\sum}\limits_{k<k_0,k\not=k(i)
}
\frac {\beta_k}{\beta_{k_0}}\,  q_k}\\
T(q_k) = q_k \quad {\rm if}\quad k\not=k_0\,.
\end{array}
$$
We set $\bar\partial_{k+1}= T^{-1}\circ \partial_{k+1} \circ T$ and 
$\bar\partial_{k}= \partial_k \circ T$.
We observe  that $\bar \partial_k=0$ on the $k$-cycles as $T$ keeps this set invariant.
We have $\bar\partial (p_j)= q_{k_0}$ and $\bar\partial(p_i) = \partial (p_i)$ if $i<j$.
Thus, $\bar\partial$ has the simple form for $1\leq i\leq j$. Arguing this way recursively, 
we get a simple complex which is $\F$-equivalent to the Morse complex.\bull

Since the equivalence relation involves upper triangular matrices only,
$\partial  (p)$ remains the class of $\De([p])$ in $H_*\bigl(f(p)-\ep;\F\bigr)$
as it is in the Morse complex. Therefore, 
when  the $\F$-equivalent complex is simple, the type of each critical point can be easily derived and this complex is the Morse-Barannikov complex. The proof of Theorem \ref{main} is completed.
\bull


\section{Bifurcations}

 Bifurcations occur in a path of functions. It follows from Thom's transversality theorems,
as it is explained by J. Cerf in the beginning of \cite{cerf},
that the space $\mathcal F$ of real smooth functions on $M$ has a natural stratification whose 
strata of
codimension $\leq 1$ are the following:
\begin{itemize}
\item[0)]  The stratum $\mathcal F_0$, an open dense set in $\mathcal F$, is formed by Morse functions whose  critical values are all simple. The next two strata are of codimension one.
\item[1)]  The stratum $\mathcal F_1$ is formed by the functions whose critical values are simple 
and  whose critical points are  all non-degenerate but one  
where the Hessian has a kernel of dimension one.
\item[2)] The stratum $\mathcal F_2$ is formed by the Morse functions whose critical values
are all simple but one which has multiplicity 2.
\item[3)] The complement $\mathcal R$  in $\mathcal F$ of the preceding strata.

\end{itemize}

A generic path has its  end points   in $\mathcal F_0$,    
 avoids $\mathcal R$ which in turn is said to be of codimension greater than 1, and crosses 
$\mathcal F_1\cup \mathcal F_2$ transversely in a finite number of points.

Since the Morse-Barannikov complex is well defined for functions in $\mathcal F_0$ only, 
it is necessary to study the bifurcation when crossing $\mathcal F_1$ and $\mathcal F_2$.
In that aim, it is convenient to introduce the Barannikov diagram and 
the Cerf diagram, which are  defined as follows.

The {\it Barannikov diagram} deals with a generic Morse function $f\in \mathcal F_0$. Let $n=\dim M$.
The vertical lines $D_k, k=0,1,\ldots, n,$ are drawn in the plane, the  $x$-coordinate of $D_k$
being $n-k$. The critical values of $f$ of index $k$ are marked on $D_k$. When the pair 
of critical points $(p, q)$ are coupled in the Barannikov complex a segment is drawn 
from $f(p)$ to $f(q)$; clearly, the slope of this segment is negative (Fig. 1). The critical values,
of homological type, in particular $\max f$ and $\min f$,  remain non-connected to another one.\\

\begin{center}
\hskip 0cm \includegraphics[scale=0.5]{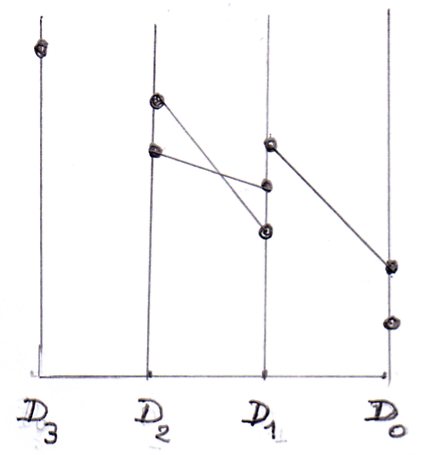}
\centerline{Figure 1}
\end{center}

The {\it Cerf diagram} deals with a generic path $\ga=\left(f_t\right)_{t\in[0,1]}$.
Its Cerf diagram is the union
 in $[0,1]\times\R$ of  $\{t\}\times f_t({\rm crit}f_t)$. It is made of finitely many smooth arcs
 transverse to the verticals  $\{t\}\times\R$, ending at cusp points  or in $\{0,1\}\times\R$,
 crossing one another transversely.
 
 \begin{rien} {\sc Bifurcation at birth times. } {\rm  This event is the crossing of the stratum 
 $\mathcal F_1$, which is co-oriented: the crossing in the positive direction corresponds to the birth of a pair of critical points of index $k$ and $k+1$ respectively. The crossing in the opposite 
 direction corresponds to the cancellation of a pair of critical points.
 
 The birth is modeled  by the following formula:
 $$f(x,y)= c+Q(y)+ x^3-(t-t_0)x,
 $$
 where $t_0$ is  the birth time, $c$ is the critical value of the birth point, 
 $(x,y)\in \R\times\R^{n-1}$
  are local coordinates at the birth point $p_0$ and  $Q$ is a non-degenerate 
  quadratic form of index $k$
   on  $\R^{n-1}$. In the Cerf diagram, there is a cusp of coordinates $(t_0, c)$. 
   For $0<t_0-t $ small, there are no critical points in the vicinity of $p_0$. 
   For $0<t-t_0$ small,
   there are two critical points $p_t, q_t$ of index $k+1$ and $k$ close to $p_0$ (Fig. 2).\\
   
   
\begin{center}
\hskip 0cm \includegraphics[scale=0.4]{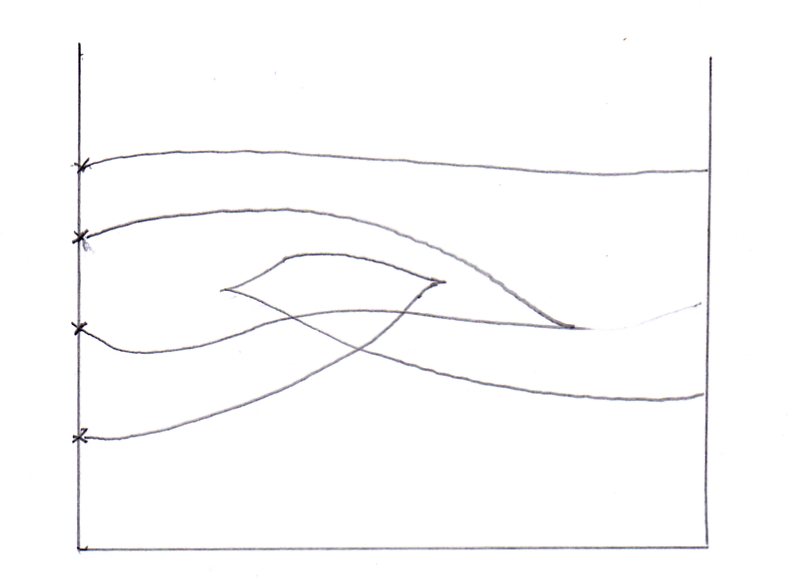}
\centerline{Figure 2}
\end{center}

   It follows from the model that $p_t$ is of upper type and $q_t$ is of lower type and this pair is coupled in the Barannikov diagram.
   The other critical points keep their type and coupling when crossing the birth time.
}
 \end{rien}

 \begin{rien} {\sc Bifurcation at double critical value time. } {\rm 
 This event is the crossing of the stratum $\mathcal F_2$, say at time $t_0$. Since the functions 
 $f_t$ are Morse for $t$ close to $t_0$, the critical points can be followed continuously on 
 $(t_0-\eta,t_0+\eta)$. Generically, the pseudo-gradient is Morse-Smale at the time $t_0$. 
 Therefore, the Morse complex remains the same on a small interval. But, as the order of the 
 critical values  is modified the Barannikov complex could change.
 
  Denote by $(p^1_t, p^2_t)$ the pair of critical points whose values
 cross at time $t_0$; say $f_t(p^1_t)<f_t(p^2_t)$ when $t<t_0$; hence, $f_t(p^1_t)>f_t(p^2_t)$ when $t>t_0$.
 The question is how the types and coupling of critical points are changing when $t$
 crosses the time $t_0$.
 
With Barannikov  we limit ourselves to the case
when $M$ is the $n$-sphere $S^n$ and  the  crossing does not involve the extremal values,
this latter question being left as an exercise to the reader.
Since on a sphere  the only critical points of  homological type are the extrema,  the crossing
deals with critical values of upper/lower type.

We shall say that there is {\it no bifurcation} if the crossing keeps all critical points with their 
initial types and coupling (remember that, near a crossing time all the functions of the considered 
path are Morse and the critical points can be followed smoothly in time due to the implicit function theorem).

 One checks  by hand that there is no bifurcation if $p^1_t$ and $p^2_t$ have distinct indices.
  Now, we are reduced to the case where the two crossing critical values have the same
   index, say $k$. 
   
  We are going to prove in the next three propositions that, in our restrictive setting, 
  there are only three
 types of bifurcations which are shown with their Barannikov diagrams before 
  and after crossing (Fig. 3
 ).\\
 \begin{rien} {\sc Notation}.
 {\rm Before stating the bifurcation propositions, it is useful to introduce some notation.
 Denote by $c_0$ the double critical value: 
 $$c_0=f_{t_0}(p_{t_0}^1)=f_{t_0}(p_{t_0}^2).
 $$
 We have two $(k-1)$-spheres $\Si^1$ and $\Si^2$ traced in the level set $f_{t_0}=c_0-\ep$
 by the respective  unstable manifolds of $p^1_{t_0}$ and $p^2_{t_0}$.
 When $t$ runs in a small interval around $t_0$, these objects move by a small isotopy:
$\Si^1_t, \Si^2_t\subset  (f_t=c_0-\ep)$.

 Finally, if $\al$ is a non-zero  $(k-1)$-homology class is the sub-level set $c_0-\ep$,
we denote by  $\la(\al)$  the critical value which is the infimum of $c$ such that 
 $\al$ vanishes in the relative homology $H_{k-1}(c_0-\ep,c; \F)$ of the pair of the 
 denoted sub-level sets.\\
}\end{rien}

 \begin{center}
\hskip 0cm \includegraphics[scale=0.7]{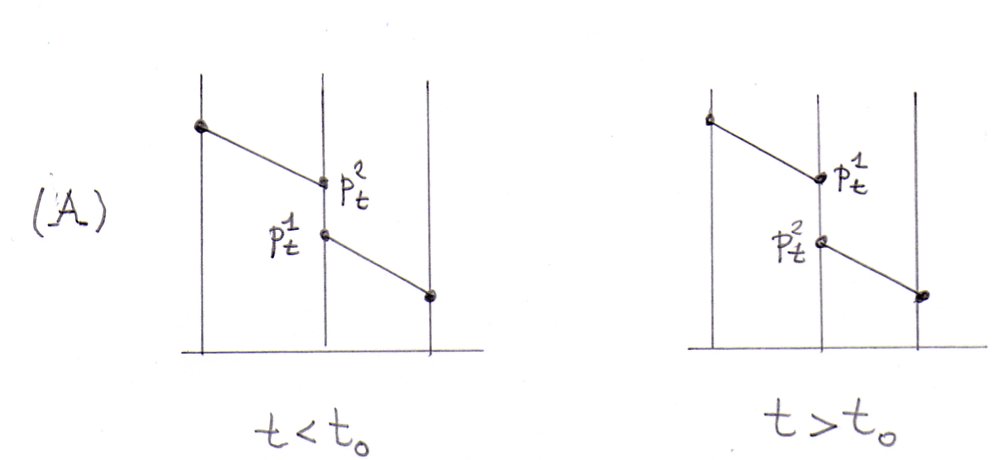}

\hskip 0cm \includegraphics[scale=0.7]{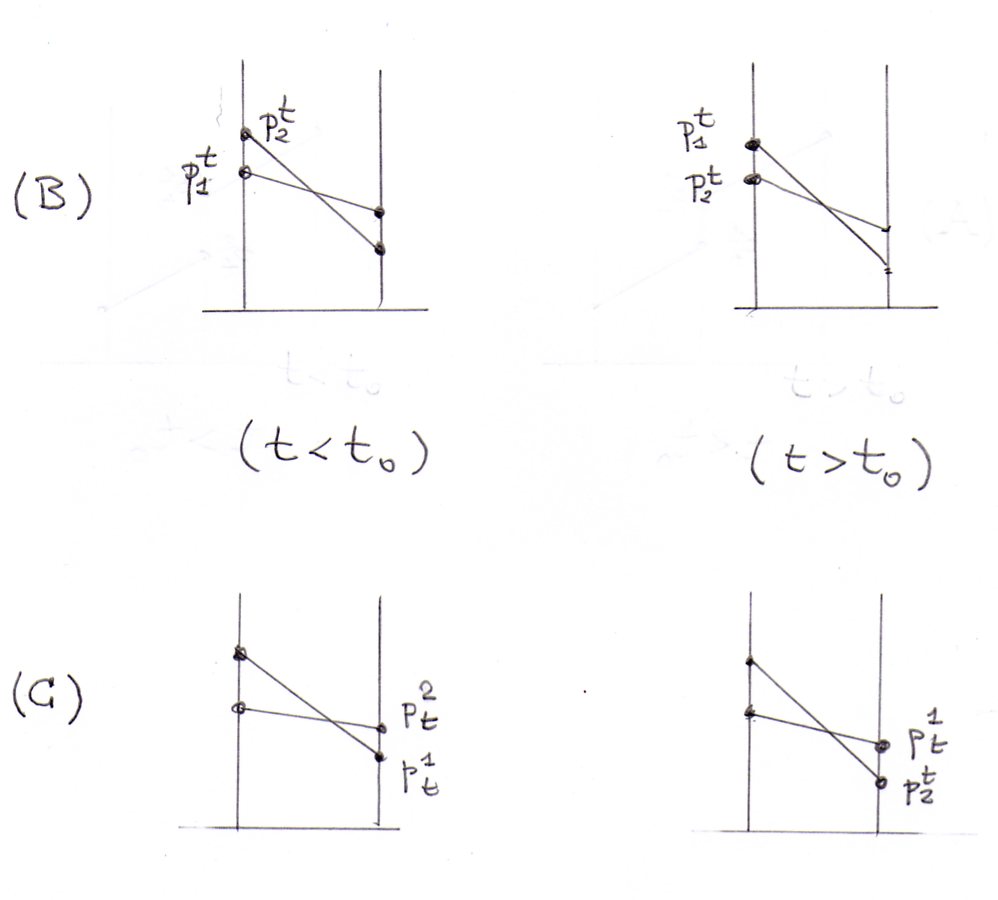}

\centerline{Figure 3}
\end{center}
 }
 \end{rien}

 \begin{prop} \label{bifA} We assume that, for $t<t_0$, $p^1_t$ and $p^2_t$ have  not the same type
 upper/lower.
 Then,  the following holds true.
 
 \nd {\rm 1)} A
 bifurcation can occur only in case {\rm (A left)}, that is: $p_t^1$ is of upper type 
 and below $p_t^2$ which is of lower type. 
 
 \nd {\rm 2)} Assuming the above necessary condition, a bifurcation occurs if and only if
  $[\Si^1]= [\Si^2]$ in $H_{k-1}(c_0;\F)$ up to a scalar.
 In that case, for $t>t_0$,
 $p^1_t$ becomes of lower type and $p^2_t$ becomes of upper type; the  new coupling is shown on {\rm (A right)}.
 \end{prop}
 
 \proof 
 1)  Assume $p^1_t$ and $p^2_t$ 
  are respectively of  lower and upper type for $t<t_0$. When $t$ is close to $t_0$, $t<t_0$,
  by assumption, $\Si^1_t$ is null-homologous
  in its sub-level set while  $\Si^2_t$ is non-homologous to 0. 
  Since these spheres change with $t$ by an  isotopy,
   this property persists up to  $t=t_0$ and still a little further, up to some $t'>t_0$. 
  Therefore the  types are unchanged after crossing and it is easy to check that the coupling is also
   unchanged. Hence, no bifurcation.

  2) Then, we assume (A left). Again, for $t<t_0$ and close to $t_0$,
  certainly $[\Si^1_t]$ is not homologous to 0 in $H_{k-1}((f_t=c_0-\ep);\F)$. The 
 unstable manifold $W^u(p^2_t)$ traces in $f_t=f_t(p^2_t)-\ep',\ \ep'>0,$ 
 a sphere which is homologous to 0 in its sub-level set if $f_t(p^2_t)-\ep'>f_t(p^1_t)$. Thus, 
  for $\Si^2_t$ in $f_t=c_0-\ep$, there are two possibilities:
   (i) $[\Si^2_t]$=0  or  (ii)  $[\Si^2_t]$ is 
  a non-zero multiple  of $[\Si^1_t]$ since it should vanish when passing in a sub-level set containing $p^1_t$. In both cases, this property persists up to $t_0$
  and up to  some $t'>t_0$.
  One checks easily that in case (i) there is no bifurcation.
  
  Consider case (ii). The point $p^2_{t'}$ is of upper type, by the very definition. Now, with
  $\ep'>0 $ small enough so that $f_{t'}(p^2_{t'})<f_{t'}(p^1_{t'})-\ep'$,
   we see that the trace of 
  $W^u(p^1_{t'})$ in the level set $f_{t'}=f_{t'}(p^1_{t'})-\ep'$ is homologous to 0 in its sub-level set.
  Therefore, $p^1_{t'}$ is of lower type and the types of the two critical points exchange. 
  One checks that  the coupling is as shown on (A right).\bull\\
  
  Assume now $p^1_t$ and $p^2_t$ are both of upper type. By the same
  homological argument as before, the type cannot change. But the coupling could  change. 
  Denote $q^1_t$ and $q^2_t$ the points associated with $p^1_t$ and $p^2_t$ respectively 
  before $t_0$.

  \begin{prop} \label{bifB} 
  In this situation, there is a bifurcation (exchange  of coupling) if and only if 
   $\la([\Si^1])= \la([\Si^2])$. In that case, when $t<t_0$,
     $f_t(q^1_t)>f_t(q^2_t)$  as shown on  {\rm (B left)}.
   \end{prop}
   
   \proof  Assume    $\la([\Si^1])< \la([\Si^2])$ (the opposite inequality is treated similarly).
   This implies that $[\Si^1]$ is 0 in the homology of the pair of sub-level sets
   $\bigl(c_0-\ep, \la([\Si^2])\bigr)$. But this vanishing holds true for every $t$ close $t_0$.
   This proves that the coupling of $p^1_t$ remains unchanged; hence, 
  no bifurcation.

What happens in case of equality? First, we look at $t<t_0$. 
Since the function  $f_t$, restricted to the sub-level set $c_0-\ep $, moves by isotopy 
when $t$ runs in $(t_0-\eta, t_0]$, one derives that 
 $\la([\Si^1_t])$ varies continuously on this interval. 
But we know $\la([\Si^1_t])=f_t(q^1_t) $; as a consequence
we have 
$$
\la([\Si^1])=\la([\Si^2])=f_{t_0}(q^1_{t_0})\,. 
$$ 
For a small $\ep'$ so that $f_t(p^2_t)-\ep'> f_t(p^1_t)$,
the trace of $W^u(p^2_t)$ in $f_t= f_t(p^2_t)-\ep'$ is homologous to 
0 in the pair of sub-level sets $(f_t(p^2_t)-\ep',f_t(q^1_t)-\ep)$;
indeed,  the class of $W^u(q^1_t)$ is null in this pair which contains $p^1_t$. Therefore,
$\la([\Si^2_t])<\la([\Si^1_t])$ and, hence, $f_t(q_t^2)<f_t(q^1_t)$ as shown on (B left).

Second, we look at $t>t_0$. Now,  $\la([\Si^2_t])$ varies continuously on $[t_0, t_0+\eta)$;
thus, this value is $f_t(q^1_t)$ due to the above equality at time $t_0$. Thus $p^2_t$ is coupled with $q^1_t$. Therefore, $p^1_t$ must 
be coupled with $q^1_t$; there is, indeed,  no other free place!\bull\\

The last case to consider is when $p_t^1$ and $p_t^2$ are both of lower type. For homological reasons, there is no change of types. Let $q_t^1$ and $q_t^2$ be the points of index $k+1$
with which they are coupled respectively when $t<t_0$.
 
At time $t_0$, denotes by $e_i, \ i=1,2$, the $k$-cell traced by the unstable manifold 
of $p^i_{t_0}$ in  the level set $f= c_0-\ep$. The map $\mu$ that we introduced in Lemma \ref{surj}
is still defined: $\mu([e_i])$ is the infimum of $c$ such that the class of $e_i$ is 0 
in the pair $(c, c_0-\ep)$ with coefficients in $\F$. Arguing similarly as in the previous proposition, one proves the following.

\begin{prop} \label{bifC} In this situation, there is a bifurcation (exchange  of coupling) if and only if 
   $\mu([e_1])= \mu([e_2])$. In that case  $f_t(q^1_t)>f_t(q^2_t)$ when $t<t_0$
    as shown on  {\rm (C left)}.
\end{prop}

 \bigskip
 
\section{The non-empty boundary case}

Following S. Barannikov, we discuss in this section the problem of extending 
 without critical points  a germ of function
given along the boundary $M$ of a compact $(n+1)$-dimensional manifold $W$, $M=\partial W$.
This setting was already considered in 1934 by Morse-van Schaack \cite{morse} where the
Morse inequalities have been formulated and proved for manifolds with non-empty 
boundary. Notice that generically a function $F:W\to\R$ is a Morse function whose critical points lie in the interior of $W$ and whose restriction to the boundary is Morse.

Actually, the problem considered by Barannikov in \cite{bara}
was more ambitious, that is, given a generic germ $\tilde f$ 
along the boundary $M$, to give a bound from below 
of the number of critical points
of any generic extension $F: W\to\R$ of $\tilde f$, as acute as possible. We focus
on $M=S^n$, the $n$-sphere, and $W=D^{n+1}$. Moreover, 
 we limit ourselves to answer the question of knowing when this bound is positive.
This problem was completely solved by S. Blank \& F. Laudenbach
\cite{blank} for $n=1$  and by 
C. Curley \cite{curley} for $n=2$.

\begin{rien} {\sc The framing.} {\rm Here we use Barannikov's terminology.
Given a Morse function $f:M\to \R$ a {\it framing}  of $f$ is the data of one vertical arrow
 at each critical value 
of $f$.
A generic germ $\tilde f: M\times [0, \ep)\to \R$ along the boundary
determines a framing according the following rule: for  the  critical
  $p$
of $f$ the arrow at $f(p)$ points up (resp. down) if $<d\tilde f(p), \vec n(p)>$ is positive
(resp. negative), where $\vec n(p)$
is a tangent vector at $p$ pointing inwards. Conversely, the framing classifies the germ $\tilde f$ up to isotopy fixing the boundary. 

This   yields some information about the non-existence of an extension 
without critical points $F:W\to\R$
of the germ $\tilde f$. For instance, if  the framing points up at the maximum, the maximum principle tells us that any extension must have at least one critical point
in the interior. Barannikov's discussion will be of course more subtle. The framing may be attached
to a Morse function on $M$,  an $\F$-equivalence  class of Morse complexes or  the  associated Morse-Barannikov complex
$C_B(f)$ as well. We will speak of  the {\it framed  Barannikov diagram.} 
We are going to look at the framed Barannikov complex and get some information in relation
 to the extension problem.}
\end{rien}

 From now on, we restrict ourselves to $M=S^n$ and $W= D^{n+1}$. A framed function
 $f:S^n\to \R$ is said to be {\it standard} if $f$ has two critical points only, the framing pointing down at the maximum and pointing up at the minimum. In this case, there is a standard 
 extension to the $(n+1)$-ball without critical points.\\

 Now, start with a generic germ $\tilde f$ along $S^n$ and assume there
 is an extension without critical points $F:D^{n+1}\to \R$. 
 Then, there is a one-parameter family of spheres 
 $S_t\subset D^{n+1},\ t\in[0,1]$, such that $S_0= M$, $S_{t'}$ lies inside $S_t$ when $t'>t$
 and the germ of $F$ along $S_1$ is {\it standard}. Set $f_t:= F\vert S_t$; it is a function 
 thought of as defined on $M=S^n$, equipped with a framing due to the knowledge
 of $S_{t'}$ for $t'=t+\ep, \ \ep>0$. The framing says the direction of moving of the spheres near a critical point $p$ of $f_t$ (up to isotopy, $F$ may be locally
 thought of  as the height function in $\R^3$).
 
 Generically, for 
 such  a family of spheres,   the map $t\mapsto f_t$
 is a generic path of functions in the sense of Section 2. 
 Therefore, there is a sequence of bifurcations starting from a given germ and ending at the 
 standard germ. 
 
 There is no bifurcation of framing. At a birth time
 the pair $(p_t,q_t)$ which has born with $f_t(p_t)>f_t(q_t)$ has the {\it standard} framing:
  the arrows are both  up or both down; such a pair will be said to be 
  a {\it standard pair}. At a crossing time $t_0$  involving the pair of critical values
  $f_t(p_t^1)<f_t(p_t^2)$, with the same indices, each critical point keeps its
  framing during the crossing. But, the rule of moving added to the rule of numbering
  implies that,  for $t<t_0$, the arrow of $p_t^1$ is up and the arrow of $p_t^2$ is down (Compare Fig. 4).
 At each generic time, we have a well-defined Morse-Barannikov framed complex
 and the bifurcations are those allowed by Propositions \ref{bifA} to \ref{bifC}. Of course, 
 in each case the bifurcation occurs only when the required homological condition is satisfied.
 If it is not, the crossing yields no bifurcation.

 \begin{center}
\hskip 0cm \includegraphics[scale=0.7]{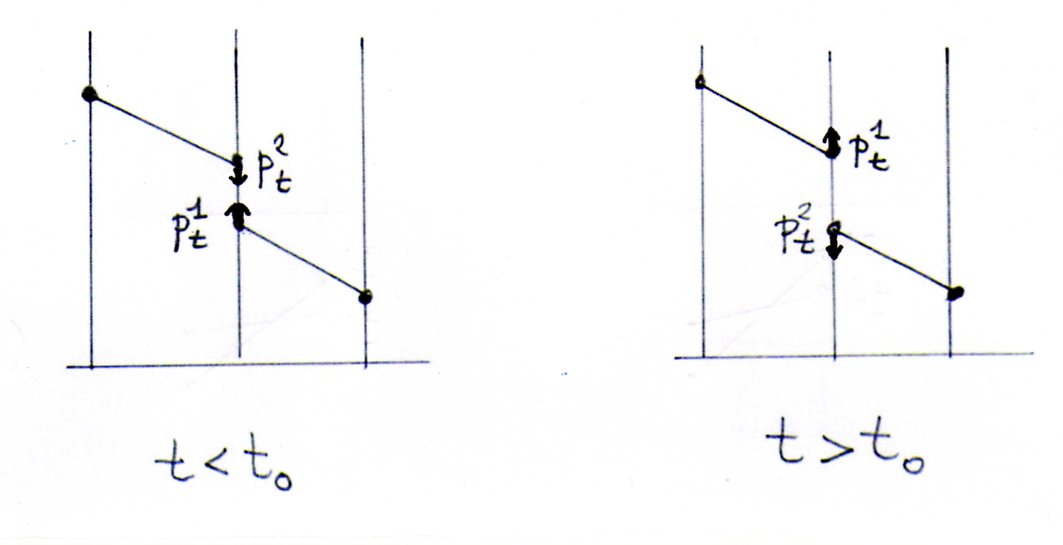}

\centerline{Figure 4}
\end{center}
 
 \medskip
 
  \begin{rien} {\sc Barannikov idea's.}
  {\rm Not any sequence of allowed bifurcations of the framed Barannikov's complex 
  is realizable by a sequence of bifurcations of framed Morse complex; at the level of framed 
  Barannikov diagram there is no longer homological condition.
  Therefore, there is a {\it formal problem} which is the following:
  given a framed Barannikov complex, does there exist a sequence of allowed bifurcations
   connecting it to the standard Barannikov complex 
   ({\it i.e.} one maximum down and one minimum up)? If there is no solution to the formal problem,
   {\it a fortiori} there is no solution to the extension problem without critical point.
   Now, the  question is whether it is possible to answer the formal problem in finite time. This 
   question is solved by the last theorem in Barannikov's article.
   }\end{rien}
   
   \begin{thm}\label{obstruction} Given a framed Barannikov complex $C_0$, if it is connected to the standard 
   Barannikov complex $C_{st}$, then $C_0$ is connected to  $C_{st}$ without any  birth.
   \end{thm}
   
   In particular,  the formal problem  reduces to a finite combinatorics.
   
   \begin{defn} Let $f$ be a framed Morse function.
   A coupled pair of critical points is said to be inverted when one of its  two arrows points up
   and the other points down. 
   
   There are two types of such inverted  pairs: 
   in type I (resp. type II) the upper point is equipped with an arrow pointing up (resp. down). 
   
   The index of a coupled pair (inverted or standard) will be the index of the upper point. 
   \end{defn}
   
   One checks on the list of bifurcations that such a pair could not disappear alone. At best,
   it is possible to shift the indices of the involved critical points (use the bifurcation of  Fig. 4  
   one pair being inverted and the other being standard in the sense of the birth bifurcation).
   But, 
   two  inverted pairs  involving of a bifurcation as shown on Fig.4  become both standard 
   and then, 
   each of them can be cancelled.
   
   As a consequence, an obvious obstruction to extending without critical points
   is the parity of the number of inverted pairs in the initial framed Morse function.
   The obstruction which follows from Theorem \ref{obstruction} is more subtle
   as shown in  the next example.
   
  \vskip 2cm  
 \begin{center}\hskip 0cm \includegraphics[scale=0.7]{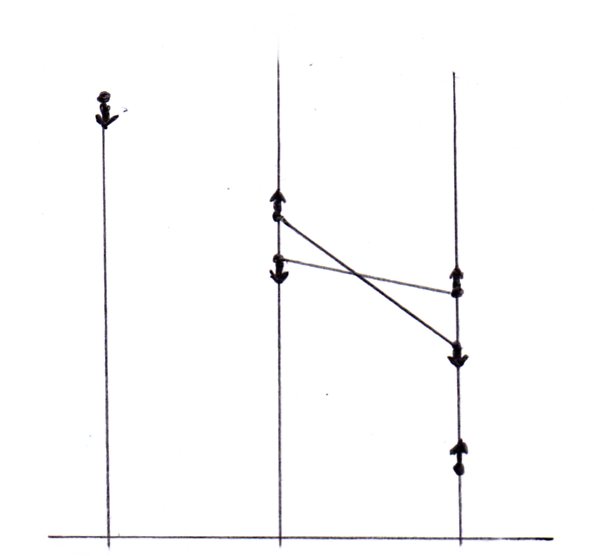}

\centerline{Figure 5}
\end{center}
\medskip

The proof of Theorem \ref{obstruction} is based on the fact that there is a very short 
list of bifurcations involving at least one inverted pair; moreover,  the role of the two types, I and II,
 are completely
different. Here is the list of these bifurcations: 
\begin{enumerate}
\item a pair of type I and a pair of type II of the same indices yield two standard pairs;

\item two pairs of type I whose indices differ by 1 yield two standard pairs whose indices differ by 1;

\item  two standard pairs having the same index $k$  yield a pair of type I and a pair of type II
both having the index $k$;

\item two standard pairs whose indices are $k$ and $k+1$ yield two pairs of type II whose indices
are still $k$ and $k+1$;

\item a pair of type I and index $k$ and a standard pair of index $k\pm 1$ yield a pair of type II
and index $k\pm 1$ and a standard pair of index $k$.
\end{enumerate}

 In particular, it is impossible to change the index of a pair of type II. For proving Theorem \ref{obstruction}, without losing generality we may assume that the initial complex $C_0$ has 
 no standard pairs. All the coupled pairs are inverted and we are facing the problem
 of canceling them, maybe with the help of introducing standard pairs (births)
 whose bifurcations could create new inverted pairs in good positions for a total cancellation. One checks that this event cannot happen.
 
 According to the previous list, for canceling one pair $A$ of type II it is required to have one pair 
 $B$
 of type I and of the same index is the right position allowing the  bifurcation (1).  If $B$ comes 
 from births followed by the bifurcation (3), then $B$ is born with another pair of type II still with the same index. So the price to pay the cancellation of $A$ with  $B$ is the appearance of $C$
 which is almost in the same position as $A$ in the Barannikov diagram, 
 up to a shift of the  heights
 in the direction of the arrows of the framing. So, nothing is gained.
The complete proof follows the same line. We refer to Barannikov's article for the details.

\end{document}